\documentclass[12pt]{article}
\begin{document}
\title{Twin Prime Sieve}
\author{H. J. Weber, Department of Physics,\\University of Virginia, 
Charlottesville, VA 22904}
\maketitle
\begin{abstract}
A sieve is constructed for ordinary twin primes of the form $6m\pm 1$ 
that are characterized by their twin rank $m.$ It does not suffer from 
the parity defect. Non-rank numbers are identified and counted using 
odd primes $p\geq 5.$ Twin- and non-ranks make up the set of positive 
integers. Regularities of non-ranks allow gathering information on 
them to obtain a Legendre-type formula for the number of twin-ranks 
at primorial arguments.             
\end{abstract}
\vspace{3ex}
\leftline{MSC: 11A41, 11N05}
\leftline{Keywords: Twin rank, twin index, non-rank numbers, sieve} 


\section{Introduction}

Our knowledge of twin primes comes mostly from sieve 
methods~\cite{hr},\cite{rm},\cite{hri},\cite{fi}. Conventional 
sieves, however sophisticated, suffer from the so-called parity 
defect: In the second member of a pair they cannot distinguish 
between a prime or product of two primes. The first of many 
real improvements of Eratosthenes' sieve was achieved by V. 
Brun~\cite{vb}. The best result for ordinary twin primes is due 
to Chen~\cite{hr},\cite{chen} proving there are infinitely many 
primes $p$ with $p+2$ either prime or a product of two primes. 

Prime numbers $p\geq 5$ are well known to be of the form~\cite{hw} 
$6m\pm 1.$ An ordinary twin prime occurs when both $6m\pm 1$ are 
prime. 

This paper is based on the original version of Ref.~\cite{ad}. 
Our goal here is to develop its mathematical foundations 
including sieve aspects and asymptotics for the twin prime 
counting function from the inclusion-exclusion principle 
applied to non-ranks (except for the remainder that is 
estimated in Ref.~\cite{hjw}). 
   
{\bf Definition~1.1.} If $6m\pm 1$ is an ordinary twin prime pair for 
some positive integer $m$, then $m$ is its {\it twin rank} and $6m$ 
its {\it twin index}. A positive integer $n$ is a {\it non-rank} if 
$6n\pm 1$ are not both prime. 

Since $2, 3$ are not of the form $6m\pm 1$ they are excluded as 
primes in the following. 

{\bf Example~1.} Twin ranks are $1,2,3,5,7,10,12,17,18,\ldots .$
Twin indices are $6,12,18,30,42,60,72,102,108,\dots .$ Non-ranks 
are $4,6,8,9,11,13,14,15,\\16,19,\ldots .$ 

In matters concerning ordinary twin primes, the natural numbers 
consist of twin- and non-ranks. Only non-ranks have sufficient  
regularity and abundance allowing us to gather enough information 
on them to draw inferences on the number of twin-ranks. Therefore, 
our main focus is on non-ranks, their symmetries and abundance.      

In Sect.~2 the twin-prime sieve is constructed based on non-ranks. 
In Sect.~3 non-ranks are identified in terms of their main properties 
and then, in Sect.~4, they are counted. In Sect.~5 twin ranks are 
isolated and then counted. Conclusions are summarized and discussed 
in Sect.~6.   

\section{Twin Ranks, Non-Ranks and Sieve}

It is our goal here to construct a twin prime sieve. To this end, 
we need the following arithmetical function. 

{\bf Definition~2.1.} Let $x$ be real. Then $N(x)$ is the integer 
nearest to $x.$ The ambiguity for $x=n+\frac{1}{2}$ with integral 
$n$ will not arise in the following.    

{\bf Lemma~2.2.} {\it Let $p\geq 5$ be prime. Then} 
\begin{eqnarray}
N(\frac{p}{6})=\{\begin{array}{ll}\frac{p-1}{6},~\rm{if}~p\equiv 
1\pmod{6};\\\frac{p+1}{6},~\rm{if}~p\equiv -1\pmod{6}.\\
\end{array}
\end{eqnarray}
{\bf Proof.} This is obvious from Def.~2.1 by substituting 
$p=6m\pm 1.~\diamond$ 

{\bf Corollary~2.3.} {\it If $p\equiv 1\pmod{6}$ is prime and 
$p-2$ is prime, then $\frac{p-1}{6}$ is a twin rank. If 
$p\equiv -1\pmod{6}$ and $p+2$ is prime, then $\frac{p+1}{6}$ 
is a twin rank.}  

{\bf Proof.} This is immediate from Def.~1.1.~$\diamond$  

{\bf Example~2.} This is the case for $p=7,13,19,31,43,61,
73,\ldots$ as well as for $p=5,11,17,29,41,59,71,\dots$ 
but not for $p=23,37,47,53,67,\dots .$  

{\bf Lemma~2.5} {\it Let $p\geq 5$ be prime. Then all 
natural numbers 
\begin{eqnarray}
k(n,p)^{\pm}=n p\pm N(\frac{p}{6})>0,~n=1,2,\ldots 
\end{eqnarray}
are non-ranks; there are $2=2^{\nu(p)}$ (single) non-rank 
progressions to $p.$  
  
(a)} If $p\equiv 1\pmod{6}$ the non-rank $k(n,p)^+$ has 
\begin{eqnarray}
6k(n,p)^+=6np+(p-1)
\label{i++}
\end{eqnarray} 
{\it sandwiched by the pair} 
\begin{eqnarray}
([6n+1]p-2,[6n+1]p),
\label{p++}
\end{eqnarray}
{\it and the non-rank $k(n,p)^-$ has}  
\begin{eqnarray}
6k(n,p)^-=6np-(p-1)
\label{p+-}
\end{eqnarray} 
{\it sandwiched by the pair} 
\begin{eqnarray}
([6n-1]p,[6n-1]p+2).
\label{p+++}
\end{eqnarray}
{\it (b) If} $p\equiv -1\pmod{6}$ the non-rank 
$k(n,p)^+$ has $6k(n,p)^+=6np+(p+1)$ sandwiched by   
\begin{eqnarray}
([6n+1]p,[6n+1]p+2);
\label{p-+}
\end{eqnarray} 
{\it and the non-rank $k(n,p)^-$ has} 
\begin{eqnarray}
6k(n,p)^-=6np-(p+1)
\label{p--}
\end{eqnarray}
{\it sandwiched by the pair}
\begin{eqnarray}
([6n-1]p-2,[6n-1]p).
\label{p---}
\end{eqnarray}
Clearly, these non-ranks are symmetrically distributed at equal 
distances $N(p/6)$ from multiples of each prime $p\geq 5.$ The 
cases for $n=0$ are the subject of Cor.~2.3 and Example~1.  
When there are more than two such non-rank progressions then 
another prime number will be involved. This issue will be 
addressed in Sect.~3. 
  
{\bf Proof.} Let $p\equiv 1\pmod{6}$ be prime and $n>0$ an integer. 
Then $k(n,p)^{\pm}=np\pm\frac{p-1}{6}$ by Lemma~2.2 and $6np\pm(p-1)$ 
are sandwiched by the two pairs in Eqs.~(\ref{p++}),(\ref{p+-}) both 
of which contain a composite number. Hence $k(n,p)^{\pm}$ are 
non-ranks.  

If $p\equiv -1\pmod{6}$ and prime, then $k(n,p)^{\pm}=np\pm\frac{p+1}
{6}$ by Lemma~2.2 and $6np\pm(p+1)$ lead to the two pairs in 
Eqs.~(\ref{p-+}),(\ref{p--}) both of which contain a composite 
number again. Hence $k(n,p)^{\pm}$ are non-ranks.~$\diamond$

The converse of Lemma~2.5 holds, i.e. {\bf non-ranks are prime 
number driven.}   
 
{\bf Lemma~2.6.} {\it If $k>0$ is a non-rank, there is a prime 
$p\geq 5$ and a non-negative integer $\kappa$ so that 
$k=k(\kappa,p)^+$ or} $k=k(\kappa,p)^-.$ 

{\bf Proof.} Let $6k+1$ be composite. Then $6k+1\neq 2^\mu 3^\nu,
~\mu,\nu\geq 1$ because then $6k+1\equiv 0\pmod{6},~$ q.e.a. If 
$\mu=0$ then $6k+1\equiv 3\pmod{6},~$ q.e.a. If $\nu=0$ then $6k+1
\equiv \pm 2\pmod{6},~$ q.e.a. If $6k+1=2^\lambda 3^\mu 5^\nu$ then 
$\lambda+\mu\leq 1,$ so $6k+1=2^\lambda 5^\nu$ and $6k+1\equiv (\pm 2)
(\pm 1)\pmod{6},~$ q.e.a. or $6k+1=3^\mu 5^\nu\equiv \pm 3\pmod{6},~$ 
q.e.a. Hence $6k+1=p\cdot K,$ where $p\geq 5$ is the smallest prime 
divisor. If $p=6m+1,$ then $K=6\kappa+1$ and  
\begin{eqnarray}
6k+1=6^2 m\kappa+6(m+\kappa)+1,~k=6m\kappa+m+\kappa=p\kappa
+\frac{p-1}{6}.  
\end{eqnarray}      
q.e.d. If $p=6m-1,$ then $K=6\kappa-1$ and  
\begin{eqnarray}
6k+1=6^2 m\kappa-6(m+\kappa)+1,~k=6m\kappa-m-\kappa=p\kappa
-\frac{p+1}{6}.
\end{eqnarray}
The case where $6k-1$ is composite is handled 
similarly.$~\diamond$ 
 
The $k(n,p)^{\pm}$ yield pairs $6k(n,p)^{\pm}\pm 1$ with one or two 
composite entries that are twin-prime analogs of multiples $np,~n>1$ 
of a prime $p$ in Eratosthenes' prime sieve~\cite{hw}. 

Concrete steps to construct a genuine prime pair sieve will 
be taken in the next section. But it is worth pointing out 
that many of the non-ranks determined using Lemma~2.5 can 
be found with the help of primes lower than $p$, but none 
of them can be found using primes greater than $p$. This 
feature is important because it ensures that, once a 
number is shown to be a twin rank by some prime up to 
a certain prime, it is not going to be shown to be a 
non-rank by any larger primes. 
 
\section{Identifying Non-Ranks} 

Here it is our goal to systematically characterize 
and identify non-ranks among natural numbers. 

{\bf Definition~3.1} Let $K>0$ be integral, 
${\cal P}_K=\{p:~~prime,~k(n,p)^+=K=k(n',p)^-\}$ 
and $p(K)=\min_{p\in {\cal P}_K}.$ Then $p(K)$ is 
called {\it parent} prime of $K.$ 

{\bf Example~3.}  
\begin{eqnarray}
p=5: k^+=6,11,16,21,\dots ;~k^-=4,9,14,19,\ldots
\end{eqnarray} 
These $k^{\pm}$ form the set ${\cal A}_5^-=\{5n\pm 
1>0:~n>0\}={\cal A}_5$ of non-ranks of parent 
prime $p=5.$ The initial non-ranks ${\cal A}_5^{(0)}=
\{4, 6\}$ give all others for $n\geq 1.$ Note that $5$ 
is the most effective non-rank generating prime number. 
If it were excluded like $3$ then many numbers, such as 
$4, 21, 56, 59, 61, 66, 74, 81, 84,\\91, 94, \ldots,$ 
would be missed as non-ranks. 

Also $(n+1)p-\frac{p-1}{6}-np-\frac{p-1}{6}=p-
\frac{p-1}{3}$ are gaps between non-ranks of the 
prime $p\equiv 1\pmod{6}.$ So are $np+\frac{p-1}{6}
-(np-\frac{p-1}{6})=\frac{p-1}{3}.$ 

For $p\equiv -1\pmod{6}$ the gaps are 
$(n+1)p-\frac{p+1}{6}-np-\frac{p+1}{6}=p-
\frac{p+1}{3}.$ And so are $\frac{p+1}{3}.$  

Thus, for $p=5$ the gaps $2,3$ in the set 
${\cal A}_5^-$ of non-ranks simply alternate.    

{\bf Proposition~3.3.} {\it The arithmetic 
progressions} $6\cdot 5n\pm 1, 6(5n+2)\pm 1, 
6(5n+3)\pm 1,~n\geq 0$ {\it contain all twin 
prime pairs except for} $3,5; 5,7.$ 

Note that the arithmetic progression $6(5n+1)+1$ 
contains $7$ of the twin $5,7$ for $n=0$ and 
infinitely many non-twin primes (by Dirichlet's 
theorem) like $37,$ while $6(5n+1)-1$ is 
composite except for $n=0.$ The set of constants 
$c\in {\cal C}_5=\{0, 2, 3\}$ of $5n+c$ in 
Prop.~3.3.   

Prop.~3.3 is the first step of the twin-prime 
sieve.

{\bf Proof.} From $\{6m\pm 1:~m=5n,5n\pm 1,5n\pm 
2,m>0\}$ we strike all pairs $6(5n+1)\pm 1,~6(5n
-1)\pm 1$ resulting from non-ranks of ${\cal 
A}_5^-.~\diamond$ 

For $p=7,$ we now subtract from the set ${\cal 
A}_7^+=\{7n\pm 1>0:~n>0\}$ of non-ranks 
the non-ranks of $p=5.$ The remaining set 
${\cal A}_7$ comprises the non-ranks to parent 
prime $p=7.$ 
 
{\bf Lemma~3.4.} {\it The set ${\cal A}_7$ of 
non-ranks to parent prime $p=7$ comprises the 
arithmetic progressions} $\{7(5n+1)+1,7(5n+2)
\pm 1,7(5n+3)\pm 1,7(5n+4)-1: n=0,1,2,\ldots\}.$

The initial non-ranks are ${\cal A}_7^{(0)}=\{8,
13,15,20,22,27\}.$ 

{\bf Proof.} The arithmetic progressions $5\cdot 
7n\pm 1\in {\cal A}_5^-,{\cal A}_7^+;$ so are  
$5(7n+1)+1=7(5n+1)-1\in {\cal A}_5^-, 
{\cal A}_7^+,~5[7(n+1)-1]-1=7(5n+4)+1\in 
{\cal A}_5^-, {\cal A}_7^+.$ Subtracting them 
from ${\cal A}_7^+,$ the arithmetic progressions 
listed in Lemma~3.4 are left over.~$\diamond$ 

Note that these four arithmetic progressions 
contain all common (double) non-ranks of the 
primes $5,~7.$ 

{\bf Proposition~3.5.} {\it The arithmetic 
progressions} $6\cdot35n\pm 1,6(35n+2)\pm 1, 
6(35n+3)\pm 1, 6(35n+5)\pm 1,6(35n+7)\pm 1,
6(35n+10)\pm 1,6(35n+12)\pm 1,6(35n+17)\pm 1,
6(35n+18)\pm 1,6(35n+23)\pm 1,6(35n+25)+1,
6(35n+28)\pm 1,6(35n+30)\pm 1,6(35n+32)\pm 1,
6(35n+33)\pm 1,~n\geq 0$ {\it contain all twin 
pairs except for} $3,5; 5,7.$ 
  
In {\bf short notation} we list the constants 
$c$ of the arithmetic progressions $35n+c$ as 
${\cal C}_7=\{c\}=\{0,2,3,5,7,10,12,17,18,23,
25,28,30,32,33\}; {\cal C}_5\subset {\cal C}_7.$   
Except for $0, 28$ all $c$ are twin ranks. 
We call the non-rank $28$ to prime $13$ an 
{\bf intruder.}

{\bf Proof.} Using Lemma~3.4, we strike from the 
arithmetic progressions of Prop.~3.3 (replacing 
$n\to 7n,7n+1,\ldots,7n+6$) all pairs resulting 
from non-ranks in ${\cal A}_7,$ which are 
$6[7(5n+1)+1]-1, 6[7(5n+2)\pm 1]-1, 6[7(5n+3)
\pm 1]\mp 1, 6[7(5n+4)-1]+1.$ This leaves the 
progressions listed above. The progressions 
$6(35n+4)\pm 1,~6(35n+6)\pm 1$ are missing 
because each of their pairs contains numbers 
divisible by $5.~\diamond$  

This is the second step of the sieve. 

From Lemma~2.5, for $p=11$ the set of non-ranks 
${\cal A}_{11}^-=\{11n\pm 2:~n=1,2,\ldots\}.$ 
First we subtract from ${\cal A}_{11}^-$ the 
non-ranks of ${\cal A}_5^-.$ Common non-ranks of 
${\cal A}_5^-$ and ${\cal A}_{11}^-$ are the 
arithmetic progressions $5(11n+2)-1=11(5n+1)-2,
5(11n+5)-1=11(5n+2)+2,5(11n+6)+1=11(5n+3)-2,
5(11n+9)+1=11(5n+4)+2.$ Note that again there 
are four common or double non-rank progressions 
to the pair of primes $5, 11$. Subtracting them 
from ${\cal A}_{11}^-$ this yields the set 
${\cal A}'_{11}=\{5\cdot 11n\pm 2,11(5n+1)+2, 
11(5n+2)-2,11(5n+3)+2,11(5n+4)-2: n\geq 0\}$ of 
arithmetic progressions. Next we subtract from 
${\cal A}'_{11}$ the common (double) non-ranks 
of ${\cal A}_7.$ Again there are four arithmetic 
progressions  
\begin{eqnarray}\nonumber
&&7(11n+2)-1=11(7n+1)+2,~n\geq 0;\\\nonumber
&&7(11n+3)-1=11(7n+2)-2,~n\geq 0;\\\nonumber
&&7(11n+8)+1=11(7n+5)+2,~n\geq 0;\\\nonumber
&&7(11n+9)+1=11(7n+6)-2,~n\geq 0.
\end{eqnarray} 
This yields ${\cal A}_{11},$ the non-ranks to 
parent prime $p=11.$ In the short notation of 
Prop.~3.5, the arithmetic progressions $2\cdot 
3(5\cdot 7\cdot 11n+c)\pm 1$ containing all 
twin primes are ${\cal C}_{11}=\{0, 2, 3, 5, 
7, 10, 12, 17, 18, 23, 25, 28, 30, 32, 33, 
37, 38, 40,\\45, 47, 52, 58, 60, 63, 65, 67, 
70, 72, 73, 77, 80, 82, 87, 88, 93, 95, 98, 
100, 102, 103,\\105, 107, 110, 115, 117, 122, 
128, 133, 135, 137, 138, 140, 142, 143, 147, 
150, 157,\\158, 165, 168, 170, 172, 173, 175, 
177, 180, 182, 187, 192, 193, 198, 203, 205, 
208,\\210, 212, 213, 215, 217, 220, 227, 228, 
235, 238, 242, 243, 245, 247, 248, 250, 252,
\\257, 263, 268, 270, 275, 278, 280, 282, 283, 
285, 287, 290, 292, 297, 298, 303, 305,\\308, 
312, 313, 315, 318, 320, 322, 325, 327, 333, 
338, 340, 345, 347, 348, 352, 353,\\355, 357, 
360, 362, 367, 368, 373, 375, 378, 380, 382\}$
except for $3,5; 5,7.$ Here $28_{13}, 37_{13}, 
60_{19}, 63_{13}, 65_{17}, 67_{13}, 73_{19}, 
\ldots$ are intruder non-ranks to the prime 
listed as subindex, while all other $c$ are 
twin ranks, except for $0$. Note that 
${\cal C}_7\subset {\cal C}_{11},$ but this 
pattern does not continue. This completes 
the 3rd step of the sieve.      

We now display characteristic properties 
of ordinary twin primes that shed light on 
the pivotal role of $N(p/6)$ and the relevance 
of non-ranks of Lemmas~2.5 and 2.6 for twin 
primes.    

{\bf Theorem~3.6.} {\it Let $p'>p\geq 5$ be 
primes such that $N(\frac{p'}{6})=N(\frac{p}
{6}).$ Then} $p'=p+2.$ 

{\bf Proof.} If $p'\equiv -1\pmod{6}$ then 
$N(\frac{p'}{6})=\frac{p'+1}{6}.$ Suppose  
$p\equiv 1\pmod{6},$ then $N(\frac{p}{6})=
\frac{p-1}{6}$ and $\frac{p'+1}{6}=
\frac{p-1}{6}.$ Hence $p'=p-2,$ q.e.a. So 
$p\equiv -1\pmod{6}$ and $N(\frac{p}{6})=   
\frac{p+1}{6}$ implies $\frac{p'+1}{6}= 
\frac{p+1}{6}.$ Hence $p'=p,$ q.e.a. So 
$p'\equiv 1\pmod{6}$ and $N(\frac{p'}{6})
=\frac{p'-1}{6}.$ Suppose $p\equiv 
1\pmod{6},$ then $N(\frac{p}{6})=
\frac{p-1}{6}.$ Since $p'>p,~N(\frac{p'}
{6})>N(\frac{p}{6}),$~q.e.a. Hence 
$p\equiv -1\pmod{6}$ and $N(\frac{p}{6})=
\frac{p+1}{6}.$ Therefore $\frac{p'-1}{6}=
\frac{p+1}{6}$ and $p'=p+2.~\diamond$ 

{\bf Corollary~3.7.} {\it Let $p\geq 5,~
p'=p+2$ be prime. Then $pp'n\pm\frac{p+1}
{6}>0$ for $n=0,1,2,\ldots$ and  
\begin{eqnarray}\nonumber
&&p(p'n+\frac{p+1}{6})+\frac{p+1}{6}=
p'(pn+\frac{p+1}{6})-\frac{p'-1}{6}>0,\\
\nonumber
&&p(p'n-\frac{p+1}{6})-\frac{p+1}{6}=
p'(pn-\frac{p+1}{6})+\frac{p'-1}{6}>0,\\
&&n=0,1,2,\ldots 
\label{nonr}
\end{eqnarray} 
are their common non-ranks.}   

Note that, again, there are four arithmetic  
progressions of common or double non-ranks. 

{\bf Proof.} Using $N(p'/6)=N(p/6)=\frac{p+1}
{6},$ Eq.~(\ref{nonr}) is readily verified;  
its lhs $\in {\cal A}_p^-$ and rhs $\in 
{\cal A}_{p+2}^+$ and $p(p+2)n\pm \frac{p+1}{6} 
\in {\cal A}_p^-,{\cal A}^+_{p+2}$.~$\diamond$ 

{\bf Example~4.} For $p'=7,~p=5$ the first 
two common non-ranks in the proof of Lemma~3.4 
are cases of Cor.~3.7. For $n=0$ small 
non-ranks are obtained.   

We now consider more systematically common  
non-ranks of pairs of primes, generalizing 
Cor.~3.7 to arbitrary prime pairs $p, p'$.  

{\bf Theorem~3.9.} {\it Let $p'>p\geq 5$ 
be primes. (i) If $p'\equiv p\equiv -1\pmod{6},$ 
then $p'=p+6l,~l\geq 1,~ N(\frac{p'}{6})=
N(\frac{p}{6})+l$ and common non-ranks of 
$p',p$ are, for} $n=0, 1,\ldots,$
\begin{eqnarray}
p(p'n+r')\pm N(\frac{p}{6})=p'(pn+r)\pm N(\frac{p'}{6})
\label{snr1}
\end{eqnarray}
{\it provided the nonnegative integers $r, r'$ solve} 
\begin{eqnarray}
(r'-r)p=l(6r\pm 1). 
\label{snr11}
\end{eqnarray}
{\it Eq.~(\ref{snr11}) with $6r\pm 1\equiv 0\pmod{p}$ 
on the rhs has a unique solution $r$ that then 
determines} $r'.$ 

{\it If $r, r'$ solve} 
\begin{eqnarray}
(r'-r)p=l(6r\mp 1)l\mp N(\frac{p+1}{6}) 
\label{snr12}
\end{eqnarray}
{\it then the common non-ranks are}
\begin{eqnarray}
p(p'n+r')\pm N(\frac{p}{6})=p'(pn+r)\mp 
N(\frac{p'}{6}).
\label{snr13}
\end{eqnarray}
{\it (ii) If $p'\equiv p\equiv 1\pmod{6},$ 
then $p'=p+6l,~l\geq 1,~N(\frac{p'}{6})=
N(\frac{p}{6})+l,$ and common non-ranks of 
$p',p$ are}
\begin{eqnarray}
p(p'n+r')\pm N(\frac{p}{6})=p'(pn+r)\pm 
N(\frac{p'}{6})
\label{snr2}
\end{eqnarray}
{\it provided $r, r'$ solve}
\begin{eqnarray}
(r'-r)p=l(6r\pm 1). 
\label{snr22}
\end{eqnarray}
{\it If $r, r'$ solve} 
\begin{eqnarray}
(r'-r)p=l(6r\mp 1)l\mp N(\frac{p-1}{6}) 
\label{snr23}
\end{eqnarray}
{\it then the common non-ranks are}
\begin{eqnarray}
p(p'n+r')\pm N(\frac{p}{6})=p'(pn+r)\mp 
N(\frac{p'}{6}).
\label{snr24}
\end{eqnarray}
{\it (iii) If $p'\equiv 1\pmod{6},~p\equiv 
-1\pmod{6}$ then $p'=p+6l+2,~l\geq 0,
~N(\frac{p'}{6})=N(\frac{p}{6})+l,$ and 
common non-ranks of $p',p$ are}
\begin{eqnarray}
p(p'n+r')\pm N(\frac{p}{6})=p'(pn+r)\pm 
N(\frac{p'}{6})
\label{snr3}
\end{eqnarray}
{\it provided}
\begin{eqnarray}
(r'-r)p=6lr+2r\pm l. 
\label{snr33}
\end{eqnarray}
{\it If $l=0$ then $r'=r=0$ and Eq.~(\ref{nonr}) 
are solutions (Cor.~3.7).} 

{\it If $r, r'$ solve} 
\begin{eqnarray}
(r'-r)p=2r(3l+1)l\mp \left(l+\frac{p+1}
{3}\right),~l\geq 1,  
\label{snr31}
\end{eqnarray}
{\it then the common non-ranks are}
\begin{eqnarray}
p(p'n+r')\pm N(\frac{p}{6})=p'(pn+r)\mp 
N(\frac{p'}{6}).
\label{snr32}
\end{eqnarray}
{\it (iv) If $p'\equiv -1\pmod{6},~p\equiv 
1\pmod{6}$ then $p'=p+6l-2,~l\geq 1,
~N(\frac{p'}{6})=N(\frac{p}{6})+l,$ and
common non-ranks of $p',p$ are}
\begin{eqnarray}
p(p'n+r')\pm N(\frac{p}{6})=p'(pn+r)\pm 
N(\frac{p'}{6})
\label{snr4}
\end{eqnarray}
{\it provided}
\begin{eqnarray}
(r'-r)p=6lr-2r\pm l. 
\label{snr44}
\end{eqnarray}
{\it If $r, r'$ solve} 
\begin{eqnarray}
(r'-r)p=2r(3l-1)\mp\left(l+\frac{p-1}{3}\right) 
\label{snr41}
\end{eqnarray}
{\it then the common non-ranks are}
\begin{eqnarray}
p(p'n+r')\pm N(\frac{p}{6})=p'(pn+r)\mp 
N(\frac{p'}{6}).
\label{snr42}
\end{eqnarray}
Note that, again, there are $4=2^{\nu(pp')}$ 
arithmetic progressions of common or double 
non-ranks to the primes $p', p$ in all cases. 
When there are more than four non-rank 
progressions then a 3rd prime will be involved. 
This case is the subject of Theor.~3.11 below.  

{\bf Proof.} By substituting $p',N(p'/6)$ 
in terms of $p,N(p/6)$ and $l,$ respectively, 
it is readily verified that Eqs.~(\ref{snr1}), 
(\ref{snr11}) are equivalent, as are 
(\ref{snr12}), (\ref{snr13}), and (\ref{snr2}), 
(\ref{snr22}), and (\ref{snr23}), (\ref{snr24}), 
and (\ref{snr3}), (\ref{snr33}), and  
(\ref{snr31}), (\ref{snr32}), and (\ref{snr4}), 
(\ref{snr44}), and (\ref{snr41}), (\ref{snr42}).  
As in (i) there is a unique solution 
$(r, r')$ in all other cases as well.~$\diamond$ 

{\bf Example~5.} For $p=5,~p'=7$ we have $l=0$ 
and Eq.~(\ref{snr33}) then gives $5(r'-r)=2r,$ 
i.e. $r'=r=0.$ Eq.~(\ref{snr3}) now gives the 
common non-ranks $35n\pm 1\in {\cal A}_5^{-}, 
{\cal A}_7^{+}.$ The common non-ranks 
$5(7n+6)-1=7(5n+4)+1$ are the solution of 
Eq.~(\ref{nonr}) with $n\to n+1.$   

For $p=5,~p'=11$ we have $l=1.$ Eq.~(\ref{snr11}) 
$5(r'-r)=6r-1$  gives the solution $r=1$ and 
$r'=r+1=2.$ Eq.~(\ref{snr1}) now gives the   
common non-ranks $5(11n+2)-1=11(5n+1)-2.$ Also  
Eq.~(\ref{snr11}), $5(r'-r)=6r+1,$ has the 
solutions $r=4,~r'=9,$ and Eq.~(\ref{snr1}) 
with plus signs displays common non-ranks of 
$5, 11,~n\geq 0.$ 

For $p'=11, p=7$ Eq.~(\ref{snr44}), $7(r'-r)=4r-1,$ 
has the solutions $r=2, r'=3$ and Eq.~(\ref{snr4}) 
with minus signs displays common non-ranks of 
$7, 11$ for $n=0, 1, 2,\ldots.$ 

For $p'=13, p=5$ Eq.~(\ref{snr33}), $5(r'-r)=8r-1,$ 
has the solutions $r=2,~r'=5,$ and Eq.~(\ref{snr3}) 
with minus signs displays common non-ranks of 
$5, 13,~n\geq 0.$ 

For $p'=13, p=7$ Eq.~(\ref{snr22}), $7(r'-r)=6r+1,$ 
has the solutions $r=1,~r'=2,$ and Eq.~(\ref{snr2}) 
with plus signs displays common non-ranks of 
$7, 13,~n\geq 0;$ and Eq.~(\ref{snr22}), $7(r'-r)
=6r-1,$ has the solutions $r=6,~r'=11,$ so 
Eq.~(\ref{snr2}) with minus signs displays 
common non-ranks of $7, 13,~n\geq 0.$ 

{\bf Theorem~3.11. (Triple non-ranks)} {\it Let 
$5\leq p<p'<p''$ (or $5\leq p<p''<p'$, or 
$5\leq p''<p<p'$) be different odd primes. Then 
each case in Theor.~3.9 of four double non-ranks 
leads to $8=2^{\nu(pp'p'')}$ triple non-ranks of 
$p, p', p''.$ At two non-ranks per prime, there 
are at most $2^3$ triple non-ranks.}  

{\bf Proof.} It is based on Theor.~3.9 and 
similar for all its cases. Let's take (i) 
and substitute $n\to p''n+\nu,~0\leq \nu<p''$ in 
Eq.~(\ref{snr1}) which, upon dropping the term 
$p''p'pn,$ yields on the lhs 
\begin{eqnarray}
pp'\nu+pr'-N(\frac{p}{6})=p''\mu\pm N(\frac{p''}
{6}).
\label{31}
\end{eqnarray} 
Since $(pp',p'')=1$ there is a unique residue $\nu$ 
modulo $p''$ so that the lhs of Eq.~(\ref{31}) is 
$\equiv \pm N(\frac{p''}{6})\pmod{p''},$ and this 
determines $\mu.$ As each sign case leads to such 
a triple non-rank solution, it is clear that there 
are $2^3$ non-ranks to $p, p', p''.~\diamond$ 
 
{\bf Example~6.} For $5,7,11$ the $2^3$ triple 
non-rank progressions are the following. Starting 
from the double non-rank equations 
\begin{eqnarray}
5(11n+2)-1=11(5n+1)-2,
\end{eqnarray}
replace $n\to 7n+\nu,$ drop $5\cdot 7\cdot 11n$ 
and set the rhs to $7\mu+1:$ 
\begin{eqnarray}
5\cdot 11\nu+5\cdot 2-1=7\mu+1.
\end{eqnarray}
Since $5\cdot 11+9=7(11-2)+1$ the solution is 
$\nu=1,~\mu=9.$ Putting back $5\cdot 7\cdot 11n$ 
we obtain the triple non-rank system 
\begin{eqnarray}\nonumber 
&&5[11(7n+1)+2]-1=7[11(5n+1)-2]+1=7[5(11n+2)-1]+1\\
&&=11[5(7n+1)+1]-2. 
\end{eqnarray}
Setting the rhs to $7\mu-1$ yields the second such  
solution  
\begin{eqnarray}
55\nu+9=7\mu-1,~\nu=3,~\mu=5^2,
\end{eqnarray}
\begin{eqnarray}
5[11(7n+3)+2]-1=7\cdot 5(11n+5)-1=11[7(5n+2)+2]-2.
\end{eqnarray}
\begin{eqnarray}
5\cdot 27-1=7\cdot 19+1=11\cdot 12+2
\end{eqnarray}
leads to
\begin{eqnarray}
5[11(7n+2)+5]-1=7[5(11n+4)-1]+1=11[5(7n+2)+2]+2; 
\end{eqnarray}
and solving the other case  
\begin{eqnarray}
55\nu+5^2-1=7\mu-1,~\nu=4,~\mu=5\cdot 7
\end{eqnarray}
leads to 
\begin{eqnarray}
5[11(7n+4)+5]-1=7\cdot 5(11n+7)-1=11[7(5n+3)+1]+2. 
\end{eqnarray}  
\begin{eqnarray}
5\cdot(11\cdot 2+6)=7\cdot 20+1=11\cdot 13-2
\end{eqnarray}
leads to  
\begin{eqnarray} 
5[11(7n+2)+6]+1=7\cdot 5(11n+4)+1=11[5(7n+2)+3]-2; 
\end{eqnarray}
and solving 
\begin{eqnarray}
5\cdot 11\nu+5\cdot 6+1=7\mu-1,~\nu=4,~\mu=36
\end{eqnarray}
leads to 
\begin{eqnarray}\nonumber
&&5[11(7n+4)+6]+1=7[5(11n+7)+1]-1=11[7(5n+3)+2]-2\\
&&=11[5(7n+4)+3]-2.
\end{eqnarray}
\begin{eqnarray}
(11\cdot 4-2)+1=7\cdot 5\cdot 6+1=11(5\cdot 4-1)+2
\end{eqnarray}
leads to 
\begin{eqnarray}
5[11(7n+4)-2]+1=7\cdot 5(11n+6)+1=11[5(7n+4)-1]+2; 
\end{eqnarray}  
and then solving
\begin{eqnarray}
5[11(\nu+1)-2]+1=7\mu-1,~\nu=5,~\mu=46
\end{eqnarray}
leads to
\begin{eqnarray}
5[11(7n+6)-2]+1=7[5(11n+9)+1]-1=11[5(7n+6)-1]-2.
\end{eqnarray}

{\bf Theorem~3.13. (Multiple non-ranks)} {\it Let 
$5\leq p_1<\cdots <p_m$ be $m$ different primes. 
Then there are $2^m$ arithmetic progressions of 
$m-$fold non-ranks to the primes} $p_1,\ldots,p_m.$ 

{\bf Proof.} This is proved by induction on $m.$   
Theors.~3.9 and 3.11 are the $m=2, 3$ cases. If 
Theor.~3.13 is true for $m$ then for any case 
$5\leq p_{m+1}<p_1<\cdots<p_m,$ or $\ldots,$ 
$5\leq p_1<\cdots<p_{m+1},$ we substitute in 
an $m-$fold non-rank equation $n\to p_{m+1}n+\nu$ 
as in the proof of Theor.~3.11, again dropping 
the $n\prod_1^{m+1} p_i$ term. Then we get   
\begin{eqnarray}\nonumber
&&p_1(p_2(\cdots(p_m\nu+r_m)+\cdots+r_2)+
N(\frac{p_1}{6})\\&&=p_{m+1}\mu\pm 
N(\frac{p_{m+1}}{6})
\label{m}
\end{eqnarray}  
with a unique residue $\nu\pmod{p_{m+1}}$ so 
that the lhs of Eq.~(\ref{m}) becomes 
$\equiv N(\frac{p_{m+1}}{6})\pmod{p_{m+1}},$ 
which then determines $\mu.$ In case the lhs of 
Eq.~(\ref{m}) has $p_1(\ldots)-N(p_1/6)$ the 
argument is the same. This yields an $(m+1)-$fold 
non-rank progression since each sign in 
Eq.~(\ref{m}) gives a solution. Hence there 
are $2^{m+1}$ such non-ranks. At two non-ranks 
per prime there are at most $2^m$ non-rank 
progressions.~$\diamond$  

{\bf Remark~3.14.}
In the multiple non-rank equations, $n\prod_1^m 
p_i$ contains the $n-$dependence, while the 
arithmetical details $r_i$ are in other 
additive terms that are independent of $n.$ 
This is the reason why the counting of 
non-ranks in the next sections will be 
independent of these arithmetical details. 

\section{Counting Non-Ranks}

If we subtract for case (i) in Theor.~3.9, say, 
the four common non-rank progressions 
corresponding to the solutions $0\leq r_i\leq 
r'_i,~0<l_i$ arranged as $0\leq r_1\leq r_2\leq 
r_3\leq r_4<p$ for definiteness, this leaves in  
${\cal A}_{p'}^-=\{p'n\pm\frac{p'+1}{6}: n\geq 0\}$ 
the following progressions $p'pn\pm \frac{p'+1}
{6},\ldots, p'(np+r_1)+\frac{p'+1}{6},\ldots, p'
(np+r_2)-\frac{p'+1}{6},\ldots, p'(np+r_3)+
\frac{p'+1}{6},\ldots, p'(np+r_4)-\frac{p'+1}
{6},\ldots, p'np\pm\frac{p'+1}{6}.$ 
 
We summarize this as follows. 

{\bf Lemma~4.1.} {\it $p'>p\geq 5$ be prime.
Removing the common non-ranks of $p',p$ 
from the set of all non-ranks of $p'$ 
leaves arithmetic progressions of the 
form $p'np+l;~n\geq 0$ where $l>0$ are 
given integers.}  

{\bf Proposition~4.2.} {\it Let $p>5$ be prime. 
Then the set of non-ranks to parent prime $p,~
{\cal A}_p,$ is made up of arithmetic 
progressions $L(p)n+a,~n\geq 0$ with $L(p)=
\prod_{5\leq p'\leq p}p,~p'$ prime and $a>0$ 
given integers.}  

{\bf Proof.} Let $p=6m\pm 1.$ We start from the 
set ${\cal A}_p^\pm=\{pn\pm N(\frac{p}{6})>0:
n=0,1,2,\ldots\}.$  Removing the non-ranks 
common to $p$ and $5$ by Lemma~4.1 leaves 
arithmetic progressions of the form $5pn+l,
~n\geq 0$ where $l>0$ are given integers. 
Continuing this process to the largest 
prime $p'<p$ leaves in ${\cal A}_p$ arithmetic 
progressions of the form $L(p)n+a,~n\geq 0$ 
with $L(p)=\prod_{5\leq p'\leq p}p'$ and $a>0$ 
a sequence of given integers independent of 
$n.~\diamond$         

The set ${\cal A}_p^{(0)}=\{a\}$ of initial 
non-ranks repeats in ${\cal A}_p$ for 
$n\geq 1$ with increment (or period) $L(p).$ 

{\bf Proposition~4.3.} {\it Let $p\geq 
p'\geq 5$ be primes and $G(p)$ the number 
of non-ranks $a\in {\cal A}_p$ over one  
period $L(p)$ corresponding 
to arithmetic progressions 
$L(p)n+a\in {\cal A}_p.$ Then} 
$G(p)=2\prod_{5\leq p'<p}(p'-2).$ 

Note that $G(p)<L(p)$ both increase 
monotonically as $p\to\infty.$ 

{\bf Proof.} In order to determine $G(p)$ 
we have to eliminate all non-ranks of 
primes $5\leq p'<p$ from ${\cal A}_p.$ 
It suffices to treat the non-ranks $a$ 
for $n=0.$ As in Lemma~3.4 we start by 
subtracting the fraction $2/5$ from the 
interval $1\leq a\leq L(p)$ of length 
$L(p),$ then $2/7$ for $p'=7$ and so on 
for all $p'<p.$ The factor of $2$ is due 
to the symmetry of non-ranks around each 
multiple of $p'$ according to Lemma~2.5.  
This leaves $p\prod_{5\leq p'<p}(p'-2)$ 
numbers $a$. The fraction $2/p$ of these 
are the non-ranks to parent prime $p,$ 
which proves Prop.~4.3.~$\diamond$   

Prop.~4.3 implies that the fraction of 
non-ranks related to a prime $p$ in the 
interval occupied by ${\cal A}_p,$  
\begin{eqnarray}
q(p)=\frac{G(p)}{L(p)}=\frac{2}{p}
\prod_{5\leq p'<p}\frac{p'-2}{p'},
\label{qp}
\end{eqnarray}
where $p'$ is prime, decreases 
monotonically as $p$ goes up.  

{\bf Definition~4.4.} Let $p\geq p'\geq 5$ be 
prime. The supergroup ${\cal S}_p=\bigcup_{
5\leq p'\leq p}\\{\cal A}_{p'}$ contains the 
sets of non-ranks corresponding to arithmetic 
non-rank progressions $a+L(p')n$ of all 
${\cal A}_{p'},~p'\leq p.$   

Thus, each supergroup ${\cal S}_p$ contains 
nested sets of non-ranks related to primes 
$5\leq p'\leq p.$

Let us now count prime numbers from $p_1=2$ 
on. 

{\bf Proposition~4.5.} {\it Let $p_j\geq 5$ 
be the $j$th prime. (i) Then the number of 
non-ranks $a\in {\cal A}_{p_i}$ 
corresponding to arithmetic progressions 
related to a prime $5\leq p_i<p_j,$}
\begin{eqnarray}
G(p_i)=\frac{L(p_j)}{L(p_i)}G(p_j)=\frac{2L(p_j)}
{p_i}\prod_{5\leq p<p_i}\frac{p-2}{p}=q(p_i)L(p_j),
\end{eqnarray}
{\it where $p$ is prime, monotonically decreases 
as $p_i$ goes up. (ii) The number of non-ranks in 
a supergroup ${\cal S}_{p_j}$ over one period 
$L(p_j)$ is}
\begin{eqnarray}
S(p_j)=L(p_j)\sum_{5\leq p\leq p_j}q(p)=L(p_j)
\left(1-\prod_{5\leq p\leq p_j}\frac{p-2}{p}
\right).\label{lsp} 
\end{eqnarray}
{\it (iii) The fraction of non-ranks of their 
arithmetic progressions in the (first) interval 
$[1,L(p_j)]$ occupied by the supergroup} 
${\cal S}_{p_j},$
\begin{eqnarray}
Q(p_j)=\frac{S(p_j)}{L(p_j)}=\sum_{5\leq p\leq p_j}
q(p)=1-\prod_{5\leq p\leq p_j}\frac{p-2}{p},
\label{qp1}
\end{eqnarray}
{\it increases monotonically as $p_j$ goes up.} 

{\bf Proof.} (i) follows from Prop.~4.3 
and Eq.~(\ref{qp}). (ii) and (iii) are 
equivalent and are proved by induction as follows, 
using Def.~4.4 in conjunction with Eq.~(\ref{qp}).  

From Eq.~(\ref{qp}) we get $q_3=2/p_3$ which is 
the case $j=3,~p_j=5$ of Eq.~(\ref{qp1}). Assuming 
Eq.~(\ref{qp1}) for $p_j,$ we add $q_{j+1}$ of 
Eq.~(\ref{qp}) and obtain
\begin{eqnarray}\nonumber
\sum_{i=3}^{j+1}q(p_i)&=&1-\prod_{i=3}^j
\frac{p_i-2}{p_i}+\frac{2}{p_{j+1}}\prod_{i=3}^j
\frac{p_i-2}{p_i}\\&=&1-\prod_{i=3}^{j+1}
\frac{p_i-2}{p_i}.
\end{eqnarray}
The extra factor $0<(p_{j+1}-2)/p_{j+1}<1$ shows 
that $q(p_j), x(p_j)$ in Eq.~(\ref{xp}) decrease 
monotonically as $p_j\to p_{j+1}$ while $Q(p_j)$ 
increases as $j\to\infty.~\diamond$

{\bf Definition~4.6.} Since $L(p)>S(p),$ there 
is a set ${\cal R}_p$ of {\it remnants} $r\in 
[1, L(p)]$ such that $r\not\in{\cal S}_p.$   

{\bf Lemma~4.7.} {\it (i) The number $R(p_j)$ of 
remnants in a supergroup, ${\cal S}_{p_j},$ is}   
\begin{eqnarray}
R(p_j)=L(p_j)-S(p_j)=L(p_j)(1-Q(p_j))=
\prod_{5\leq p\leq p_j}(p-2)=\frac{1}{2}G(p_{j+1}).
\label{rps}
\end{eqnarray} 

{\it (ii) The fraction of remnants in} 
${\cal S}_{p_j},$   
\begin{eqnarray}
x(p_j)&=&\frac{R(p_j)}{L(p_j)}=1-Q(p_j)=
\prod_{5\leq p\leq p_j}\frac{p-2}{p},
\label{xp}
\end{eqnarray}
{\it where $p$ is prime, decreases monotonically 
as} $p_j\to\infty.$ 

{\bf Proof.} (i) follows from Def.~4.6 in 
conjunction with Eq.~(\ref{lsp}) and (ii) 
from Eq.~(\ref{rps}). Eq.~(\ref{rps}) follows 
from Eq.~(\ref{qp1}).~$\diamond$
       
\section{Remnants and Twin Ranks}

When all primes $5\leq p\leq p_j$ and appropriate 
nonnegative integers $n$ are used in Lemma~2.5 
one will find all non-ranks $k<M(j+1)\equiv 
(p_{j+1}^2-1)/6.$ By subtracting these non-ranks 
from the set of positive integers $N\leq M(j+1)$ 
all and only twin ranks $t<M(j+1)$ are left among 
the remnants, i.e. twin primes with index 
$T<p_{j+1}^2-1.$ If a non-rank $k$ is left then 
$6k\pm 1$ must have prime divisors that are 
$>p_j$ according to Lemma~2.5, which is 
impossible. 

{\bf Definition~5.1.} All $t<M(j+1)=(p_{j+1}^2-1)/6$ 
in a remnant ${\cal R}_{p_j}$ of a supergroup 
${\cal S}_{p_j}$ are twin ranks. These twin ranks 
are called {\it front twin ranks}. They are included 
in $R_0.$  

{\bf Example~7.} For $p_{18}=61,~p_{19}=67$ and 
$M(19)=748$ we get the set of remnants $m=1, 2, 3, 5, 7, 
10, 12, 17, 18, 23, 25, 30, 32, 33, 38, 40, 45, 47, 52, 58,\\
70, 72, 77,87, 95, 100, 103, 107, 110, 135, 137, 138, 143, 
147, 170, 172, 175, 177,\\ 182, 192, 205, 213, 215, 217, 
220, 238, 242, 247, 248, 268, 270, 278, 283, 287, 298,\\
312, 313,322, 325, 333, 338, 347, 348, 352, 355, 357, 
373, 378, 385, 390, 397, 425,\\432, 443, 448, 452, 455, 
465, 467, 495, 500, 520, 528, 542, 543, 550, 555, 560, 
562,\\565, 577, 578, 588, 590, 593, 597, 612, 628, 637, 
642, 653, 655, 667, 670, 675, 682,\\688, 693, 703, 705, 
707, 710, 712, 723, 737, 747,$ which are all twin ranks, 
i.e., $6m\pm 1$ are prime pairs.    

{\bf Proposition~5.3.} {\it Let $p_j$ be the $j$th  
prime number and $L(p_j)n+a_i^{(j)}$ be the arithmetic 
progressions that contain the non-ranks 
$a_i^{(j)}\in{\cal A}_{p_j}$ to parent prime $p_j.$ 

(i) Let $6[L(p_j)n+c_i^{(j)}]\pm 1$ be the arithmetic 
progressions that contain the ordinary twin primes 
with $c_i^{(j)}\in{\cal C}_{p_j}.$ If $c_i^{(j)}$ is 
a twin rank or intruder non-rank to a prime $p>p_{j+1},$ 
then $c_i^{(j)}\in {\cal C}_{p_{j+1}},$ if it is a 
non-rank to $p_{j+1}$ then} $c_i^{(j)}\not\in 
{\cal C}_{p_{j+1}}$.  

{\it (ii) The set of constants $c_i^{(j+1)}$ of 
arithmetic progressions containing the twin ranks 
from the next supergroup ${\cal S}_{p_{j+1}}$ is} 
\begin{eqnarray}\nonumber
&&{\cal C}_{p_{j+1}}=\{6[L(p_j)(p_{j+1}n+l)+c_i^{(j)}]
\pm 1\} \\&&-\{6[L(p_j)(p_{j+1}n+l')+a_{i'}^{(j)}]
\pm 1\}. 
\label{rec}
\end{eqnarray}
{\it If there are positive integers $0\leq l,l'<p_{j
+1},$ a non-rank $a_{i'}^{(j)}\in {\cal A}_{p_j}$ and 
a constant $c_i^{(j)}\in {\cal C}_{p_j}$ satisfying} 
\begin{eqnarray}
L(p_j)l+c_i^{(j)}=L(p_j)l'+a_{i'}^{(j)}, 
\end{eqnarray}
{\it then} 
\begin{eqnarray}
L(p_j)l+c_i^{(j)}\not\in {\cal C}_{p_{j+1}},
\end{eqnarray}
{\it else}
\begin{eqnarray} 
c_{i,l}^{(j+1)}=L(p_j)l+c_i^{(j)}\in {\cal 
C}_{p_{j+1}}. 
\end{eqnarray}   
Prop.~5.3 is the inductive step completing the 
practical sieve construction for ordinary twin 
primes. Props.~3.3, 3.5 and Lemma~3.4 are the 
initial steps.  

{\bf Proof.} Replacing $n\to p_{j+1}n$ in (i)  
we obtain $6[L(p_j)p_{j+1}n+c_i^{(j)}]\pm 1|_{n=0}
\in {\cal C}_{p_{j+1}},$ if $c_i^{(j)}$ is a twin 
rank or non-rank to a prime $p>p_{j+1}$. If it is 
non-rank to $p_{j+1}$ then $c_i^{(j)}\not\in 
{\cal C}_{p_{j+1}},$ which proves (i).  

Replacing in (ii) $n\to p_{j+1}n+l,~l=1, 2,\ldots, 
p_{j+1}-1$ and subtracting the resulting sets 
from each other, we obtain (ii).$~\diamond$  

For $p_3=5,$ Prop.~5.3 is Prop.~3.3, for $p_4=7$ 
it is Prop.~3.5, for $p_5=11$ the sequence of $c$ is 
listed after Prop.~3.5. Clearly, at the start of the 
$c$ for $p_4=7$ the previous values for $p_3=5$ are 
repeated, and this is also the case for $p_5=11;$ 
but this pattern does not continue, as shown in 
Prop.~5.3. 

{\bf Example~8.} For prime number $p_3=5$ the 
$(c, n)$ are given in Prop.~3.3, viz. $(1,0), (1,1), 
(1,2), (2,0), (2,1), (2,2), (3,0), (3,3),\ldots$ 
corresponding to the twin primes $6\pm 1, 30\pm 1, 
60\pm 1, 12\pm 1, 42\pm 1, 72\pm 1, 18\pm 1, 
108\pm 1, \ldots,$ respectively. 

For $p_4=7$ the additional $(c, n)$ are given in 
Prop.~3.5, viz. $(10,0), (17,0),\\(18,0), (23,0), 
(25,0), (32,0), (33,0),\ldots$ corresponding to 
the twin primes $60\pm 1, 102\pm 1, 108\pm 1, 
138\pm 1, 150\pm 1, \pm 1, 180\pm 1, 192\pm 1, 
198\pm 1, \ldots,$ respectively.   

For $p_5=11$ the additional $(c, n)$ are given 
in short notation after Prop.~3.5, viz. $(37,4), 
(38,0), (40,0), \ldots$ corresponding to the twin 
primes $9462\pm 1, 228\pm 1, 240\pm 1,\ldots,$ 
respectively.   

For $p_{18}=61$ the twin ranks are listed in 
Example~7.  
 
Twin ranks are located among the remnants 
${\cal R}_p$ for any prime $p\geq 5$. The main 
goal in this Sect.~5 is to establish the 
inclusion-exclusion principle for non-ranks 
and use it to derive the twin prime version 
of Legendre's formula for $\pi(x)-\pi(\sqrt{x})$ 
extracted from Eratosthenes' sieve~\cite{hr},\cite{fi}. 
The prime $p_j$ here plays the role of the 
variable $\sqrt{x}$ there, and the front twin 
ranks here correspond to the primes $p<\sqrt{x}$ 
left over after striking out their multiples 
there.   

{\bf Theorem~5.5.} {\it Let $R_0$ be the number 
of remnants of the supergroup ${\cal S}_{p_j},$ 
where $p_j$ is the $j$th prime number,  
$M(j+1)=[p_{j+1}^2-1]/6$ and $x=L(p_j)-M(j+1).$ 
Then the number of twin ranks within the 
remnants of the supergroup ${\cal S}_{p_j}$ is 
given by} 
\begin{eqnarray}
\pi_2(6x+1)=R_0+\sum_{p_j<n\le x,n|L_j(x)}
\mu(n)2^{\nu(n)}\bigg[\frac{x}{n}\bigg]+O(1),
\label{tr} 
\end{eqnarray}
{\it where $L_j(x)=\prod_{p_j<n\le x}p,$ 
and $O(1)$ accounts for the less than perfect 
cancellation at low values of $x$ of $R_0$ 
and the sum in Eq.~(\ref{tr}).} This 
cancellation is the subject of Theors.~5.7, 
5.8 and proved there, up to a remainder 
estimated in Ref.~\cite{hjw}.  

The upper limit in the sum~(\ref{tr}) is $x$ 
because $[x/n]=0$ for $n>x.$ Here $L(p_j)
=\prod_{5\leq p\leq p_j}p,$ and $R_0=\prod_{
5\leq p\leq p_j}(p-2)$ with $p$ prime 
includes the front twin ranks, and $n$ runs 
through all products of primes $p_j<p\le x$. 

The argument of the twin-prime counting 
function $\pi_2$ is $6x+1$ because, if 
$x$ is the last twin rank of the interval 
$[1, L(p_j)],$ then $6x\pm 1$ are the 
corresponding twin primes. The twin pair 
$3, 5$ has no twin rank and is not part 
of the remnants. 

The formula~(\ref{tr}) may be regarded as the 
twin-prime analog of Legendre's application of 
Erathostenes' sieve to the prime counting 
function $\pi(x)$ in terms of M\"obius' 
arithmetic function which was subsequently 
improved by many others~\cite{hri}.  

{\bf Proof.} According to Prop.~4.5 the 
supergroup ${\cal S}_{p_j}$ has $S(p_j)=L(p_j)\\
\cdot\left(1-\prod_{5\leq p\leq p_j}\frac{p-2}
{p}\right)$ non-ranks. Subtracting these from 
the interval $[1, L(p_j)]$ that the supergroup 
occupies gives $R_0=\prod_{5\leq p\leq p_j}(p-2)$ 
for the number of remnants which include twin 
ranks and non-ranks to primes $p_j<p\leq x .$ 
The latter are 
\begin{eqnarray}
M(j+1)<pn\pm N(\frac{p}{6})\leq L(p_j),~
M(j+1)=\frac{p_{j+1}^2-1}{6}, 
\end{eqnarray} 
or 
\begin{eqnarray}
0<n\leq \frac{L(p_j)-M(j+1)}{p}, 
\end{eqnarray}
which have to be subtracted from the remnants 
to leave just twin ranks. Correcting for 
double counting of common non-ranks to two 
primes using Theor.~3.9, of triple non-ranks 
using Theor.~3.11 and multiple non-ranks 
using Theor.~3.13 we obtain 
\begin{eqnarray}\nonumber
\pi_2(6x+1)&=&R_0-2\sum_{p_j<p\leq x,n|L_j(x)}
\bigg[\frac{x}{p}\bigg]\\&+&4\sum_{p_j<p<p'\leq x}
\bigg[\frac{x}{pp'}\bigg]\mp\cdots+O(1), 
\label{tr1}
\end{eqnarray} 
where $[x]$ is the integer part of $x$ as 
usual and where $L_j(x)=\prod_{p_j<p\le x}p.$ 
Equation~(\ref{tr1}) is equivalent to 
Eq.~(\ref{tr}).~$\diamond$ 

{\bf Definition~5.6} Replacing the floor function $[x]$ 
in Eq.~(\ref{tr1}) by its argument minus fractional part, 
$[x]=x-\{x\},$ we call 
\begin{eqnarray}\nonumber  
R_M&=&R_0+\sum_{p_j<n\le x,n|L_j(x)}\mu(n)2^{\nu(n)}
\frac{x}{n},\\R_E&=&-\sum_{p_j<n\le x,n|L_j(x)}
\mu(n)2^{\nu(n)}\bigg\{\frac{x}{n}\bigg\}+O(1)
\label{tr2}
\end{eqnarray}
the main and error terms of Eq.~(\ref{tr}) or 
Eq.~(\ref{tr1}). 

{\bf Theorem~5.7.} {\it The main term $R_M$ in 
Eq.~(\ref{tr}) satisfies} 
\begin{eqnarray}\nonumber
&&R_M=L(p_j)\prod_{5\leq p\leq x}\left(1-
\frac{2}{p}\right)\\&&+M(j+1)[1-\prod_{p_j<p
\le x}(1-\frac{2}{p})].
\label{tr3}
\end{eqnarray}

{\it Proof.} Expanding the product  
\begin{eqnarray}
R_0=L(p_j)\prod_{5\leq p\leq p_j}(1-\frac{2}{p})
\end{eqnarray}
and combining corresponding sums of $R_M$ in 
Eq.~(\ref{tr2})
\begin{eqnarray}
-\sum_{5\leq p\leq p_j}\frac{1}{p}-\sum_{
p_j<p\le x}\frac{1}{p}=-\sum_{5\leq p\leq x}
\frac{1}{p},\ldots 
\end{eqnarray} 
just shifts the upper limit $p_j$ in the product 
$\prod_{p\le p_j}(1-2/p)$ to $x$ so that we obtain 
Eq.~(\ref{tr3}).~$\diamond$ 

{\bf Theorem~5.8.} {\it The main term $R_M$ 
obeys the asymptotic law} 
\begin{eqnarray}
R_M\sim \frac{c_2 e^{-2\gamma}6x}
{\log^2(6x+1)},  
\end{eqnarray}
{\it as the $j$th prime $p_j\sim \log x\to\infty$ 
where $c_2$ is the twin prime constant.} 

{\it Proof.} The ratio of the second term in 
Eq.~(\ref{tr3}) to the leading first term is 
of order,  
\begin{eqnarray}
\frac{p_{j+1}^2\log^2 x}{x}\to 0,~p_{j+1}=O(\log 
x). 
\end{eqnarray}   

Using the prime-number theorem~\cite{hw},\cite{rm} 
we have 
\begin{eqnarray}
\log L(p_j)=\sum_{5\le p\leq p_j}\log p=p_j+R(p_j)
=\log x+O(\frac{\log^2 x}{x})
\end{eqnarray}
for all sufficiently large primes $p_j,$ 
$R(p_j)$ is the remainder in the prime number 
theorem and the error term comes from $M(j+1).$  

Using Mertens' asymptotic formula~\cite{rm} for 
$x\to\infty:$ 
\begin{eqnarray}
\prod_{p=2}^x\left(1-\frac{1}{p}\right)\sim 
\frac{e^{-\gamma}}{\log x},~\gamma\approx 0.5772,  
\end{eqnarray} 
and the twin-prime constant 
\begin{eqnarray}
c_2=\frac{\prod_{p>2}(1-\frac{2}{p})}{\prod_{p>2}
(1-\frac{1}{p})^2}
\end{eqnarray}
where $p$ runs through primes, we obtain  
\begin{eqnarray}\nonumber
&&\prod_{p>2}^x (1-\frac{2}{p})=\frac{
\prod_{p>2}^x(1-\frac{2}{p})}{\prod_{p>2}^x
(1-\frac{1}{p})^2}\prod_{p>2}^x(1-\frac{1}{p})^2
\sim c_2\prod_{p>2}^x(1-\frac{1}{p})^2\\&&\sim
\frac{4c_2 e^{-2\gamma}}{\log^2 x},~x\to\infty  
\end{eqnarray} 
and thus 
\begin{eqnarray}
R_M\sim \frac{c_2 e^{-2\gamma} 6x}
{\log^2(6x+1)},~\log x\to\infty  
\end{eqnarray}
from Theorem~5.7.~$\diamond$

\section{Summary and Discussion} 

The twin prime sieve constructed here differs 
from other more general sieves that are applied 
to twin primes among many other problems in that 
it is conceptually designed for ordinary prime 
twins. It is specific rather than general and 
has no precursor in sieve theory. Hence much 
of the work is devoted to developing its 
concepts into sieve tools. Accurate 
counting of non-rank sets require the 
infinite primorial set $\{6L(p_j)=\prod_{
p\leq p_j}p\}$.  
 
The twin primes are not directly sieved, 
rather twin ranks $m$ are with $6m\pm 1$ 
both prime. All other natural numbers 
are non-ranks. These are much more 
numerous and orderly than twin ranks. 
Surprisingly, their order is governed 
by all primes $p\geq 5.$ In contrast to 
other sieves, primes serve to organize  
and classify non-ranks in arithmetic 
progressions with equal distances 
(periods) that are primes ($\geq 5$) 
or products of them. 

The coefficient $c_2 e^{-2\gamma}\approx 
0.416213$ in the asymptotic law of the 
main term $R_M$ of $\pi_2(6x+1)$ in 
Theor.~5.8 is a little less than a third of 
the Hardy-Littlewood constant $2c_2\approx 
1.320320.$ In Ref.~\cite{hjw} the remainder 
$R_E$ is shown to be at most of the order of 
the main term divided by any positive power 
of $\log x$. To put the deviation from the 
Hardy-Littlewood law in perspective, our 
``minimal'' asymptotic law holds only near 
primorials. In the large gaps between those 
special arguments there is room for other 
asymptotic laws. 

Finally, needless to say, the genuine  
sieve has no consequences for other twin 
primes or the Goldbach problem~\cite{hw}.

 
\end{document}